\documentclass{article}

\usepackage{amsmath}
\usepackage{amssymb}
\usepackage{vmargin}
\usepackage{wasysym}
\usepackage[active]{srcltx}
\usepackage[stable]{footmisc}
\usepackage{caption}
\usepackage{pdfpages}
\usepackage[matrix, arrow,all,cmtip,color]{xy}
\usepackage{url}
\usepackage{epigraph}

\newcommand{\bi}{\begin{itemize}}
\newcommand{\ei}{\end{itemize}}
\newcommand{\bd}{\begin{description}}
\newcommand{\ed}{\end{description}}

\def\lra{\longrightarrow}
\def\ra{\rightarrow}
\def\rtt{\,\rightthreetimes\,}
\def\xra{\xrightarrow}

\def\Z{\Bbb Z}
\begin{document}

\title{Expressing the statement of the Feit-Thompson theorem 
with diagrams in the category of finite groups}
\author{Misha Gavrilovich {\tt 
mi\!\!\!ishap\!\!\!p@sd\!\!\!df.org}
\thanks{
 http://mishap.sdf.org/mints-lifting-property-as-negation\newline
Higher School of Economics, 
Soyza Pechatnikov str., 16, St.Petersburg, Russia. \newline 
St. Petersburg Institute for Economics and Mathematics of 
the Russian Academy of Sciences.
}
}

\maketitle

\setlength{\epigraphwidth}{0.5\textwidth}
\epigraph{\normalsize There's no point in being grown up if you can't be childish sometimes
}{ {\emph{ 
}}
}

\begin{abstract}
\normalsize
We reformulate the statement of the Feit-Thompson theorem 
in terms of  diagrams in the category of finite groups, 
namely iterations  of the Quillen lifting property
with respect to particular morphisms. 
\end{abstract}

\Large

\section{Introduction. Structure of the Paper}

We reformulate several 
 notions in finite group theory 
in terms of diagrams in the category of finite groups,
namely being solvable, nilpotent, $p$-group and prime-to-$p$ group, 
abelian, perfect, subnormal subgroup, injective and surjective homomorphism 
(see Fig.~1 and Fig.~2). 
These notions are enough to formulate the Feit-Thompson theorem.

These properties are obtained by iterating the same diagram chasing ``trick''
in the category of groups, 
often starting from a single morphism-counterexample; the ``trick'' 
is known as  the Quillen lifting property and was introduced 
 by Quillen [Qui] to axiomatise algebraic topology 
in terms of categories. 

 Our motivation was to formulate 
part of finite group theory in a form amenable 
to automated theorem proving while remaining human readable;
[G1] tried to do the same thing for the basics of general topology. 



In more detail, we reformulate these notions using the following four operations 
$P^{\rtt l}$, $P^{\rtt r}$, $ P_{0\ra *}$, $P_{*\ra 0}$  on the classes (properties) $P$ of  morphisms in a category. 
 {\em Left, resp. right, Quillen negation} $P^{\rtt l}$, resp. $P^{\rtt r}$, of a property $P$ is 
the class of all morphisms $p^l$, resp. $p^r$,  such that
$p^l\rtt p$, resp. $p\rtt p^r$,  for all $p\in P$; see Fig.~1 for the definition of the Quillen lifting property $\rtt$: 
 $$ P^{\rtt l}:=\{\, p^l\ :\ p^l\,\rtt\,p\,\ \text{for each }p\in P\,\}$$
$$P^{\rtt r}:=\{\, p^r\ :\ p\, \rtt p^r\,\ \text{for each }p\in P\,\}$$
The classes $P^{\rtt l}$, $P^{\rtt r}$ are subcategories which contain no morphisms from $P$ but isomorphisms; hence the terminology.

Classes $P_{0\ra *}$, $P_{*\ra 0}$  denote the subclass consisting of morphisms from/to the terminal object $0$, 
i.e. the trivial group.

Classes $(P_{0\ra *})^{\rtt l}_{0\ra *}=(P_{0\ra *})^l_{0\ra *}$, 
and $(P_{0\ra *})^{\rtt r}_{0\ra *}=(P_{0\ra *})^r_{0\ra *}$
can described as classes of objects admitting no non-trivial map from/to 
the objects corresponding to $P_{0\ra *}$.

We also say that a morphism $f$ is {\em left, resp.~right, Quillen unlike $P$} or {\em $P$-unlike}
iff $f\in P^{\rtt l}$, resp.~$f \in P^{\rtt r}$, and that $f$ is 
 {\em left, resp.~right, exemplified by $P$} 
iff $f\in (P^{\rtt r})^{\rtt l}$, resp.~$f \in (P^{\rtt l})^{\rtt r}$.

In this notation our main observations are that a finite group $G$ is soluble
iff  $0\lra G$ lies in the class 
$(((\{ {\Bbb F}_2\lra {\Bbb Z}_2 \})^{\rtt rl})_{0\ra *})^{\rtt r}$
or, equivalently, 
right exemplified by taking the commutator subgroup, i.e. lies in 
 $\{ [S,S]\lra S\,:\, S\text{ an arbitrary group}\,\}^{\rtt lr}$.

A finite group $G$ is nilpotent iff the diagonal map $G\xra\Delta G\times G$, $g\mapsto (g,g)$ 
is exemplified by morphisms from the trivial group, i.e. iff 
the diagonal map lies in $(0\lra *)^{\rtt lr}$ where $(0\lra *)$ denotes the class of morphisms from the trivial group $0$. 

A finite group $G$ is perfect iff $0\lra G$ is exemplified by  $(\{ {\Bbb F}_2\lra {\Bbb Z}_2 \})$
or, equivalently, left unlike taking the commutator subgroup, 
i.e.~$0\lra G$ lies in
 $(\{ {\Bbb F}_2\lra {\Bbb Z}_2 \})^{\rtt rl}$
or, equivalently, in the class $\{ [S,S]\lra S\,:\, S\text{ an arbitrary group}\,\}^{\rtt l}$.

A group $G$ is a prime-to-$p$-group, resp. a $p$-group, iff $0\lra G$ lies in the class 
$(0\lra\Z/\!p\Z)^{\rtt r}$, resp. $(0\lra\Z/\!p\Z)^{\rtt rr}$; in words,  
it is right $(0\lra\Z/\!p\Z)$-unlike, resp. right unlike a $(0\lra\Z/\!p\Z)$-unlike morphism.
Classes $(0\lra \Z)^{\rtt r}$ and $(\Z\lra 0)^{\rtt r}$ are the classes of surjective and injective morphisms,
i.e.~a map is surjective, resp.~injective, iff it is $(0\lra \Z)$-unlike, resp.~$(\Z\lra 0)$-unlike. 

Little attempt has been made to go beyond these examples. 
Hence open questions remain: are there other interesting examples of lifting properties in the category of (finite) groups?
Can a complete group-theoretic argument be reformulated in terms of diagram chasing, say 
the classification of CA-groups or $pq$-groups, or 
elementary properties of subgroup series; can category theory notation
be used to make expositions easier to read? 
Can these reformulations be used in automatic theorem proving? 
Is there a decidable fragment of (finite) group theory
 based on the Quillen lifting property and, more generally,
 diagram chasing, cf.~[GLS]?

Examples of lifting properties outside of group theory may be found in 
a  short note [DMG] in  {\em The De Morgan Gazette} and [G1].

This approach was motivated in part by hopes to express some statements
of finite group theory in terms of diagrams in the category of finite groups,
and then use automated diagram chasing to construct short formal proofs. 
Particularly important for our motivation was the fact that in the few examples
of lifting properties we did find, it was the lifting property with respect 
to a simple counterexample in some intuitive sense, and that the lifting properties
are closely related to the usual definitions.

\def\rrt#1#2#3#4#5#6{\xymatrix{ {#1} \ar[r]^{} \ar@{->}[d]_{#2} & {#4} \ar[d]^{#5} \\ {#3}  \ar[r] \ar@{-->}[ur]^{}& {#6} }}
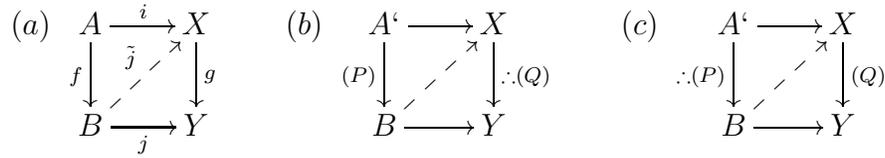
\begin{figure}
\begin{center}
\large
$(a)\  \xymatrix{ A \ar[r]^{i} \ar@{->}[d]_f & X \ar[d]^g \\ B \ar[r]_-{j} \ar@{-->}[ur]^{{\tilde j}}& Y }$
$\ \ \ \ \ \ (b)\ \xymatrix{ A \ar[r] \ar@{->}[d]_{(P)}` & X \ar[d]^{\therefore (Q)} \\ B \ar[r] \ar@{-->}[ur]& Y }$
$\ \ \ \ \ \ (c)\ \xymatrix{ A \ar[r] \ar@{->}[d]_{\therefore(P)}` & X \ar[d]^{ (Q)} \\ B \ar[r] \ar@{-->}[ur]& Y }$

\end{center}
\caption{\label{fig1}\normalsize
(a) 
The definition of a lifting property $f\rtt g$: for each $i:A\lra X$ and $j:B\lra Y$
making the square commutative, i.e.~$f\circ j=i\circ g$, there is a diagonal arrow $\tilde j:B\lra X$ making the total diagram
$A\xra f B\xra {\tilde j} X\xra g Y, A\xra i X, B\xra j Y$ commutative, i.e.~$f\circ \tilde j=i$ and $\tilde j\circ g=j$.\newline
We say that $f$ lifts wrt $g$, $f$ left-lifts wrt $g$, or $g$ right-lifts wrt $f$.
\ (b) Right Quillen negation. 
The diagram defines a property $Q$ of morphisms in terms of a property $P$; a morphism has property (label) $Q$ iff it right-lifts 
wrt any morphism with property $P$, i.e. 
$Q=\{\, p^r\ :\ p\, \rtt p^r\,\ \text{for each }p\in P\,\}$
\ (c) Left Quillen negation. 
The diagram defines a property $P$ of morphisms in terms of a property $Q$; a morphism has property (label) $P$ iff it left-lifts 
wrt any morphism with property $Q$, i.e. 
 $ P=\{\, p^l\ :\ p^l\,\rtt\,q\,\ \text{for each }q\in Q\,\}$
}
\end{figure}

\def\nto{\not\to}
\def\ZpZ{{\Bbb Z}/\!p{\Bbb Z}}
\section{Expressing properties of finite groups as diagrams.} 
Figure~2 lists the diagrams representing the properties of groups. Below we describe the same diagrams, and make a couple of remarks about the diagram chasing and the notion of an inner automorphism and Sylow theory.

 There is no non-trivial homomorphism from a group $F$ to $G$, write $F\nto G$, iff  
$$0\lra F\,\rtt 0\lra G\text{ or equivalently }F\lra 0 \rtt G\lra 0.$$
A group $A$ is {\em Abelian} iff 
\[ \left<a,b\right> \,\lra\,\left<a,b:ab=ba\right>  \rtt\,\, A\lra 0
\]
where $\left<a,b\right> \,\lra\,\left<a,b:ab=ba\right>$ is the  abelianisation morphism sending the free group into the Abelian free group on two generators; 
a group $G$ is {\em perfect}, $G=[G,G]$, iff $G\nto A$ for any Abelian group $A$, i.e. 
\[ \left<a,b\right> \,\lra\,\left<a,b:ab=ba\right>  \rtt\,\, A\lra 0\ \implies\ G\lra 0 \rtt A\lra 0\]
equivalently, for an arbitrary homomorphism $g$, 
\[ \left<a,b\right> \,\lra\,\left<a,b:ab=ba\right>  \rtt\,\,g \ \ \implies\ G\lra 0 \rtt\, \,g\,\]
Yet another reformulation is  that, for each group $S$, 
$$
0\lra G \,\rtt\, [S,S]\lra S .$$
In the category of finite or algebraic groups,
a group $H$ is {\em soluble} iff  $G\nto H$ for each perfect group $G$, 
i.e. 
$$0\lra G\,\rtt 0\lra H\text{ or equivalently }G\lra 0 \rtt H\lra 0.$$
Alternatively, a group $H$ is {\em soluble} iff for every homomorphism $f$ it holds
$$ f \,\rtt\, [G,G]\lra G \text{ for each group }G\  \implies\ f \,\rtt 0\lra\, H 
.$$
A prime number $p$ does not divide the number  elements of a finite group $G$ 
iff  $G$ has no element of order $p$, i.e. no element $x\in G$ such that $x^p=1_G$ yet $x^1\neq 1_G,...,x^{p-1}\neq 1_G$,
equivalently $\ZpZ\nto G$, i.e.
$$0\lra \ZpZ\,\rtt 0\lra G\text{ or equivalently }\ZpZ\lra 0 \rtt G\lra 0.$$
A finite group $G$ is a $p$-group, i.e. the number of its elements is a power of a prime number $p$, iff in the category of finite groups 
$$0\lra \ZpZ \rtt 0\lra H \implies 0\lra H\rtt 0\lra G.$$

A group $H$ is the normal closure of the image of $N$, 
i.e.~no proper normal subgroup of $H$ contains the image of $N$, iff for an arbitrary group $G$
$$N\lra H \,\rtt 0 \lra G.$$
A group $D$ is a subnormal subgroup of a finite group $G$ iff  
$$ N\lra H \,\rtt 0 \lra B \text{ for each group }B\ \implies  N\lra H\,\rtt\,D \lra G
$$
i.e. $D\lra G$ right-lifts wrt any map $N\lra H$ such that $H$ is the normal closure of the image of $N$; the lifting property implies that $D\lra G$ is injective. Recall that $D$ is a subnormal subgroup of a finite group $G$ iff there is a finite series of subgroups
$$D=G_0 \vartriangleleft G_1  \vartriangleleft \ldots  \vartriangleleft G_n =G$$ 
such that $G_i$ is normal in $G_{i+1}$, $i=0,\ldots,n-1$.
This is probably the only claim which requires a proof. First notice that if $D$ is normal in $G$ then 
the lifting property holds.
Given a square corresponding to $N\lra H \,\rtt\, D\lra G$, the preimage of $D$ in $H$ 
is a normal subgroup of $H$ containing the image of $N$, hence the preimage of $D$ contains $H$ and the lifting property holds. 
The lifting property is closed under composition, hence it holds for subnormal subgroups as well.
Now assume $D$ is not subnormal in $G$. As $G$ is finite, there is a minimal subnormal subgroup $D'>D$ of $G$.
By construction 
no proper normal subgroup of $D'$ contains  $D$
but the lifting property $D\lra D' \,\rtt\, D\lra G$ fails.

Finally, a finite group $G$ is nilpotent iff the diagonal group $G$ is subnormal in $G\times G$ [Nilp], i.e.~iff
the diagonal map $G\xra \Delta G\times G$, $g\mapsto (g,g)$
 right-lifts wrt any $N\lra H$ such that $H$ is the normal closure of the image of $N$,  
$$ N\lra H \,\rtt 0 \lra B \text{ for each group }B\ \implies  N\lra H\,\rtt\, G\xra \Delta G\times G
.$$

Sylow theorem implies in a  finite group $G$, each $p$-subgroup is contained in a maximal one $Syl_p(G)$, 
$Ord\, G/Ord\, Syl_p(G)$ is prime to $p$,  
and the maximal $p$-subgroups are conjugated by an inner automorphism.

It is not clear how to express this in a satisfactory manner in terms of category theory (diagram chasing).
Perhaps something along the following lines:
(in the category of finite groups) 
each arrow $0\lra G$ decomposes as 
$$0\xra{(p\text{-group})} Syl_p(G) \xra{(\text{prime to }p)} G$$ 
uniquely up to conjugation. 
Here
label $A\xra{(p\text{-group})} B $ may mean something like $Ord\, B/Ord\, Im\, A$ is a power of $p$,
and label $B \xra{(\text{prime to }p)} C$ may mean something like $Ord\, C/Ord\, Im\, B$ is prime to $p$.

We remark that the notion of 
an inner automorphism can be reformulated in a diagram chasing manner. 
An inner automorphism $g\longmapsto a g a^{-1}$ of a group $G$ extends to an automorphism $h\longmapsto \iota(a)h\iota(a)^{-1}$
of a group $H$ for any embedding $\iota:G\lra H$. [Inn, Sch] show this is a characterisation: 
an automorphism $\sigma :G\lra G$ is inner iff it extends to an automorphism of $H$ for any embedding $\iota:G\longrightarrow H$. 
See [Inn] and references therein for several more similar reformulations.

The Feit-Thompson theorem can be expressed as a combination of lifting properties:
the theorem says says that each (finite) group of odd order is soluble, i.e.
for each perfect finite group $G$ and each finite group $H$, 
$$
0\lra {\Bbb Z}/2{\Bbb Z} \rtt 0\lra H \implies 0\lra G \rtt 0\lra H.$$

Note that all these examples but the last one have a flavour of negation---a notion 
being defined by the lifting property
with respect to the simplest counterexample.

\def\rrtt#1#2#3#4#5#6#7{\xymatrix{ {#1} \ar[r]^{} \ar@{->}[d]_{#2} & {#4} \ar[d]_{#5}^{#6} \\ {#3}  \ar[r] \ar@{-->}[ur]^{}& {#7} }}
\begin{figure}
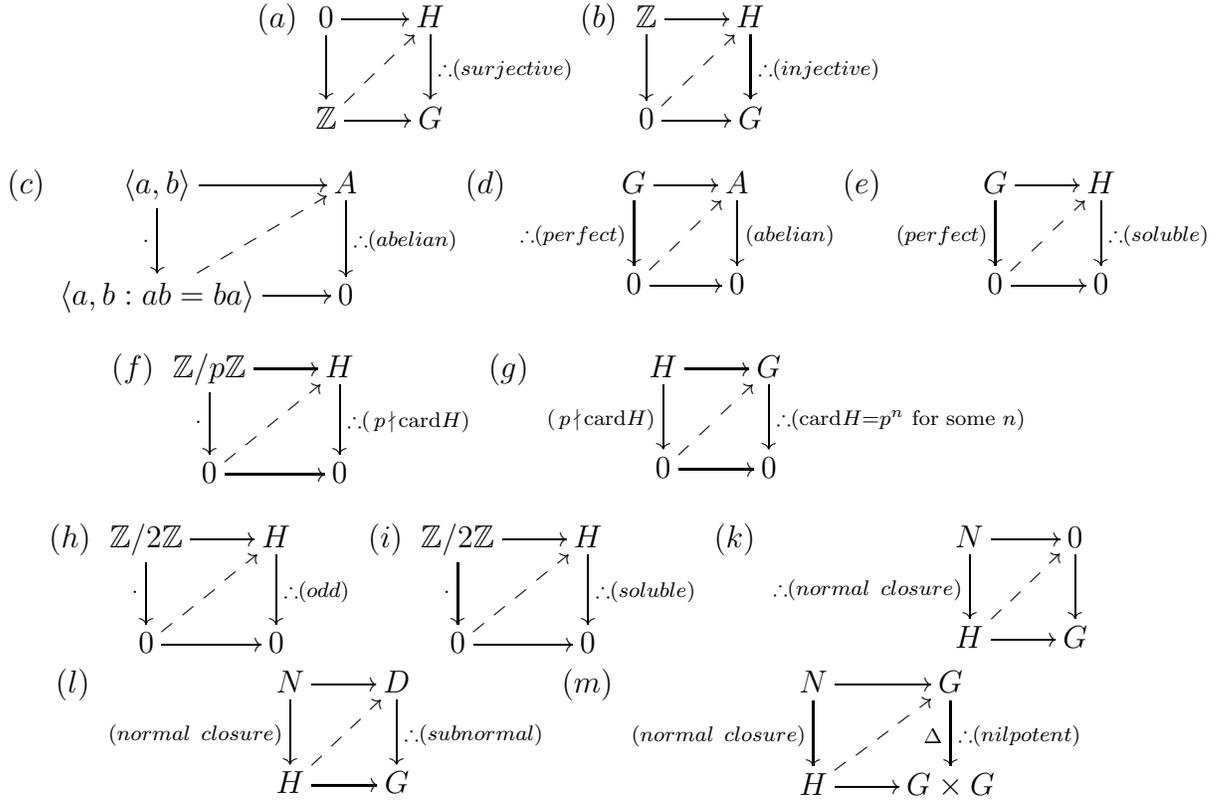

\begin{center}
\large
$(a)\  \rrt  {0} {} {\Bbb Z}   H {\therefore(surjective)} G $
$(b)\  \rrt   {\Bbb Z} {} 0  H {\therefore(injective)} G $
\vskip3mm
$(c)\  \rrt  { \left<a,b\right>} {.} {\left<a,b:ab=ba\right>}   A {\therefore(abelian)} 0 $
$(d)\  \rrt {G} {\therefore(perfect)} {0}  A {(abelian)} 0 $
$(e)\  \rrt G {(perfect)} {0} {H} {\therefore(soluble)} {0}$\ 
\vskip3mm
$(f)\  \rrt {{\Bbb Z}/p{\Bbb Z}} {.} {0} {H} {\therefore(\,p\,\nmid\,\textrm{card}H)} {0}$
$(g)\ \rrt {H} {(\,p\,\nmid\,\textrm{card}H)} {0} {G} {\therefore(\textrm{card}H=p^n\text{ for some }n)} {0}$\\  
\vskip3mm
$(h)\  \rrt {{\Bbb Z}/2{\Bbb Z}} {.} {0} {H} {\therefore(odd)} {0}$
$(i)\ \rrt {{\Bbb Z}/2{\Bbb Z}} {.} {0} {H} {\therefore(soluble)} {0}$\
$(k)\ \rrt  {N} {\therefore(normal\ closure)} {H} {0} { } {G}$
$(l)\ \rrtt   {N} {(normal\ closure)} {H} {D} {} { \therefore(subnormal)  } {G}$
$(m)\ \rrtt   {N} {(normal\ closure)} {H} {G} {\Delta} { \therefore(nilpotent)  } {G\times G}$
\end{center}
\caption{\label{fig5}\normalsize
Lifting properties/Quillen negations. Dots $\therefore$ indicate free variables. 
Recall these diagrams represent rules in a diagram chasing calculation
and ``$\therefore(label)$" reads as: 
given a (valid) diagram, add label $(label)$ to the corresponding arrow. 
A diagram is valid
iff for every commutative square of solid arrows
with properties indicated by labels, 
there is a diagonal (dashed) arrow making the total diagram commutative. 
A single dot indicates that the morphism is a constant.\newline
(a) a homomorphism $H\lra G$ is surjective, i.e.~for each $g\in G$ there is $h\in H$ sent to $g$\newline  
(b) a homomorphism $H\lra G$ is injective, i.e.~the kernel of $H\lra G$ is the trivial group\newline
(c) a group is abelian iff each morphism from the free group of two generators 
 factors through its abelianisation ${\Bbb Z}\times {\Bbb Z}$.\newline
 (d)  a group $G$ is perfect, $G=[G,G]$, iff it admits no non-trivial homomorphism to an abelian group\newline
 (e) a finite group is soluble  iff it admits no non-trivial homomorphism from a perfect  group; 
 more generally, this is true in any category of groups with a good enough dimension theory.\newline
(f) by Cauchy's theorem, a prime $p$ divides the number of elements of a finite group $G$ 
iff the group contains an element $e,e^p=1, e\neq1$ of order $p$\newline
(f) a group has order $p^n$ for some $n$ iff iff the group contains no element $e,e^l=1, e\neq1$ of order $l$ prime to $p$\newline
(h)  by Cauchy's theorem, a finite group has an odd number of elements iff it contains no involution $e,e^2=1, e\neq1$\newline
(i) The Feit-Thompson theorem says that each group of odd order is soluble,~i.e.~it says that this diagram chasing 
 rule is valid in the category of finite groups. Note that it is not a definition of the label unlike the other  
 lifting properties.\newline
(k) a group $H$ is the normal closure of the image of $N$ iff $N\lra H \,\rtt 0 \lra G$ for an arbitrary group $G$\newline
(l) $D\lra G$ is injective and the subgroup $D$ is a subnormal subgroup 
of a finite group $G$ iff  $D \lra G$
 right-lifts wrt any map $N\lra H$ such that $H$ is the normal closure of the image of $N$\newline
(m) a group $G$ is nilpotent iff the diagonal map $G\xra \Delta G\times G$, $g\mapsto (g,g)$
 right-lifts wrt any inclusion of a subnormal subgroup $N\lra H$
}\end{figure}

\section*{Acknowledgments and historical remarks}  
It seems embarrassing to thank anyone for ideas so trivial, and
we do that in the form of historical remarks....
This work is a continuation of [DMG]; early history is given there. 

Examples here were motivated by a discussion with S.Kryzhevich.
I thank Paul Schupp for pointing out the characterisation of inner
automorphisms of [Sch]. I thank M.Bays, K.Pimenov, V.Sosnilo and S.Synchuk for proofreading,
and several students for encouraging and helpful discussions. 

Special thanks are due to M.Bays for helpful discussions.

 I  wish to express my deep thanks to Grigori Mints, to whose memory this paper is dedicated \dots

\end{document}